\documentclass{amsart}

\usepackage{lscape,graphicx}
\usepackage{graphicx}

\usepackage{tikz}
\usetikzlibrary{arrows}
\tikzstyle{block}=[draw opacity=.7,line width=1.4cm]

\theoremstyle{definition}

\theoremstyle{remark}

\numberwithin{equation}{section}

\begin{document}

\title{A closed-form expression for $\zeta(2n+1)$ reveals a self-recursive function}
%\\ Cancer history, biology and systems biology

\author{{Michael A. Idowu}\\
\\
{SIMBIOS Centre, \\
University of Abertay, \\
Dundee DD1 1HG, UK\\
m.idowu@abertay.ac.uk}}

\maketitle

\begin{abstract}
Euler discovered a formula for expressing the value of the Riemann zeta function for all even positive integer arguments.  A closed-form expression for the Riemann zeta function for all odd integer arguments, based on the values of the Dirichlet beta function, euler numbers and $\pi$, reveals a new evidence about the self-recursive nature of Riemann zeta function at odd integers. We demonstrate for the first time that the Riemann zeta function at odd integers always produces a recurrence relation that is self-recursive.\\

\end{abstract}

%\tableofcontents

\section{ Introduction}
Nearly all number theorists have sought for a closed-form expression for $\zeta(2n+1)$; n being a positive integer number. Our investigation into this open problem has uncovered a self-recursive function intrinsic in $\zeta(2n+1)$. We develop and present a new method for deriving a closed-form expression for $\zeta(2n+1)$ which is based on the values of the Dirichlet beta function at $2n+1$ and Euler numbers at $2n$. This new method may be regarded as a general formula for finding a closed-form self-recursive expression for $\zeta(2n+1)$. 

Here we demonstrate how to obtain a closed-form expression for $\zeta(2n+1)$ with examples on $\zeta(3)$, $\zeta(5)$, $\zeta(7)$, $\zeta(9)$, and $\zeta(11)$. Each result produced is an exact representation, but the resultant expression is self-recursive.

\subsection{Riemann zeta and Dirichlet beta functions}
The Dirichlet beta function \cite{Nis11} is defined as 
\begin{equation}\label{def}
\beta(s)=\sum_{k=1}^{\infty} {(-1)^{k-1}\over{(2k-1)^s}} = 1-{1 \over {3^s}}+{1 \over {5^s}}-{1 \over {7^s}}+\dots
\end{equation}
\ref{def} implies 
\begin{equation}
 \beta(s)=\sum_{k=1}^{\infty} {(-1)^{k-1}\over{(2k-1)^s}} = (1+{1 \over {3^s}}+{1 \over {5^s}}+{1 \over {7^s}}+\dots) - 2({1 \over {3^s}}+{1 \over {7^s}}+{1 \over {11^s}}+\dots);
\end{equation}
$Re(s)>0$. Therefore,
\begin{equation}\label{bet}
\beta(s)= {{2^s-1}\over{2^s}}\zeta(s)-2\sum_{k=1}^{\infty} {1\over{(4k-1)^s}} 
\end{equation}
\begin{equation}%\label{beta2}
\beta(s)= {{2^s-1}\over{2^s}}\zeta(s)-{{2}\over{2^s.2^s}}\sum_{k=0}^{\infty} {1\over{(k+{3\over4})^s}}.
\end{equation}
The following equations are then obtained according to \cite{Ido12b}
\begin{equation}\label{closedform_}
\zeta(s)={{2^s}\over{2^s-1}}\beta(s)+(-1)^s{{2}\over{2^s(2^s-1)}}{1 \over{\Gamma(s)}}{{ \psi^{(s-1)}({3\over4})}};
\end{equation}
\begin{equation}\label{closedform}
\zeta(2s+1)={{2^{2s+1}}\over{2^{2s+1}-1}}\beta(2s+1)-{{2}\over{2^{2s+1}(2^{2s+1}-1)}}{1 \over{\Gamma(2s+1)}}{{ \psi^{(2s)}({3\over4})}},
\end{equation}
where the polygamma function $\psi^{(s-1)}(x)$ is defined as 
\begin{equation}
\psi^{(s-1)}(x) = {d^{s-1}\over{dx^{s-1}}}\psi(x) = {d^{s}\over{dx^{s}}}\ln \Gamma(x).
\end{equation}

\section{Closed-form expressions for $\zeta(3)$ , $\zeta(5)$, $\zeta(7)$,
 $\zeta(9)$, $\zeta(11)$}
Equation \ref{closedform} may be used to find the closed-form expressions for $\zeta(3)$, $\zeta(5)$, $\zeta(7)$, $\zeta(9)$, and $\zeta(11)$.

\subsection{A closed-form expression for the Riemann zeta of 3}
From \ref{closedform} we derive
 
\[
\zeta(3)={{2^3}\over{2^3-1}}\beta(3)-{{2}\over{2^3(2^3-1)}}{1 \over{\Gamma(3)}}{{ \psi^{''}({3\over4})}}={8\over7}\beta(3)-{{2}\over{8(7)}}{1 \over{2!}}{{(2{\bf(1)}\pi^3 - 2({\bf{28}})\zeta(3))}}
\]
\begin{equation}
\rightarrow \zeta(3)={8\over7}\beta(3)-{{1}\over{56}}{{(2\pi^3 - 56\zeta(3))}}
\end{equation}

\subsection{A closed-form expression for the Riemann zeta of 5} Similarly:
\[
\zeta(5)={{2^5}\over{2^5-1}}\beta(5)-{{2}\over{2^5(2^5-1)}}{1 \over{\Gamma(5)}}{{ \psi^{''''}({3\over4})}}={32\over31}\beta(5)-{{2}\over{32(31)}}{1 \over{4!}}{{(8{\bf(5)}\pi^5 - 8{\bf{(1488)}}\zeta(5))}}
\]
\begin{equation}
\rightarrow \zeta(5)={32\over31}\beta(5)-{{2}\over{32(31)}}{1 \over{4!}}{{(40\pi^5 - 11904\zeta(5))}}
\end{equation}

\subsection{A closed-form expression for the Riemann zeta of 7}
\[
\zeta(7)={{2^7}\over{2^7-1}}\beta(7)-{{2}\over{2^7(2^7-1)}}{1 \over{\Gamma(7)}}{{ \psi^{(6)}({3\over4})}}={128\over 127}\beta(5)-{{2}\over{128(127)}}{1 \over{6!}}{{(32{\bf(61)}\pi^7 - 32{\bf{(182880)}}\zeta(7))}}
\]
\begin{equation}
\rightarrow \zeta(7)={128\over127}\beta(7)-{{2}\over{128(127)}}{1 \over{6!}}{{(1952\pi^7 - 5852160\zeta(7))}}
\end{equation}

\subsection{A closed-form expression for the Riemann zeta of 9}
\[
\zeta(9)={{2^9}\over{2^9-1}}\beta(9)-{{2}\over{2^9(2^9-1)}}{1 \over{\Gamma(9)}}{{ \psi^{(8)}({3\over4})}}={512\over 511}\beta(9)-{{2}\over{512(511)}}{1 \over{8!}}{{(128{\bf(1385)}\pi^9 - 128{\bf{(41207040)}}\zeta(9))}}
\]
\begin{equation}
\rightarrow \zeta(9)={512\over511}\beta(9)-{{2}\over{512(511)}}{1 \over{8!}}{{(177280\pi^9 - 5274501120\zeta(9))}}
\end{equation}

\subsection{A closed-form expression for the Riemann zeta of 11}
\[
\zeta(11)={{2^{11}}\over{2^{11}-1}}\beta(11)-{{2}\over{2^{11}(2^{11}-1)}}{1 \over{\Gamma(11)}}{{ \psi^{(10)}({3\over4})}}\]
\[\rightarrow \zeta(11)= {2048\over 2047}\beta(11)-{{2}\over{2048(2047)}}{1 \over{10!}}{{(512{\bf(50521)}\pi^{11} - 512{\bf{(14856307200)}}\zeta(11))}}
\]
\begin{equation}
\rightarrow \zeta(11)={2048\over2047}\beta(11)-{{2}\over{2048(2047)}}{1 \over{10!}}{{(25866752\pi^{11} - 7606429286400\zeta(11))}}
\end{equation}

\section{The general (closed-form expression) formula for $\zeta(2s+1)$}
The results obtained in the previous sections indicate the following general formula for obtaining representing $\zeta(2s+1)$: 

\begin{equation}\label{zetaOdd}
 \zeta(2s+1)={2^{2s+1}\over(2^{2s+1}-1)}\beta(2s+1)-2{{{\biggr(2^{2s-1}{ \bf{\mid E_{2s} \mid}} \pi^{2s+1} - {{{2^{2s-1}{\bf{2(2^{2s+1}-1)}}}} {\bf{\Gamma(2s+1)}}}\zeta(2s+1)\biggr)}} \over{{{2^{2s+1}(2^{2s+1}-1)}} \Gamma(2s+1)}}
\end{equation}
; s is an integer. The modulus $\mid E_{2s} \mid$ is the absolute value of an even-indexed Euler number $E_{2s}$. The implication of \ref{zetaOdd} is 

\begin{equation}
{{{2^{2s+1}(2^{2s+1})}} \Gamma(2s+1)}\beta(2s+1)={{{2^{2s}{ \bf{\mid E_{2s} \mid}} \pi^{2s+1}}} }
\end{equation}
as expected. Hence, 
\begin{equation}
 \zeta(2s+1)={2^{2s+1}\over(2^{2s+1}-1)}\beta(2s+1)-2{{{\biggr(2^{2s-1}{ \bf{ {E_{2s}\over i} }} (\pi {\bf{i}})^{2s+1} - {{{2^{2s-1}{\bf{2(2^{2s+1}-1)}}}} {\bf{\Gamma(2s+1)}}}\zeta(2s+1)\biggr)}} \over{{{2^{2s+1}(2^{2s+1}-1)}} \Gamma(2s+1)}}
\end{equation}
$\rightarrow$ 
\begin{equation}
 \zeta(2s+1)={2^{2s+1}\over(2^{2s+1}-1)}\beta(2s+1)+{{{\biggr(2^{2s-1}{ \bf{ E_{2s} }} (\pi i)^{2s+1} - {{{2^{2s-1}{\bf{2(2^{2s+1}-1)}}}} {\bf{\Gamma(2s+1)}}}\zeta(2s+1){\bf{i}}\biggr)}} \over{{{2^{2s+1}(2^{2s+1}-1)}} \Gamma(2s+1)}}2{\bf{i}}
\end{equation}
$\rightarrow$

\begin{equation}\label{zetaOdd2}
 \zeta(2s+1)={2^{2s+1}\over(2^{2s+1}-1)}\beta(2s+1)+\biggr[{{{{ \bf{ E_{2s} }} (\pi i)^{2s+1}}} \over{{{2(2^{2s+1}-1)}} \Gamma(2s+1)}} - \zeta(2s+1)i\biggr]i
\end{equation}

\begin{equation}
{(2^{2s+1}-1)\over{2^{2s+1}}} \zeta(2s+1)=\beta(2s+1)+{(2^{2s+1}-1)\over{2^{2s+1}}}\biggr[{{{{ \bf{ E_{2s} }} (\pi i)^{2s+1}}} \over{2(2^{2s+1}-1) \Gamma(2s+1)}} - \zeta(2s+1)i\biggr]i
\end{equation}

\begin{equation}\label{last}
{(2^{2s+1}-1)\over{2^{2s+1}}} \zeta(2s+1)=\beta(2s+1)+{(2^{2s+1}-1)\over{2^{2s+1}}}\biggr[\zeta(2s+1)+{{{{ \bf{ E_{2s} }} (\pi i)^{2s+1}}} \over{2 (2^{2s+1}-1)\Gamma(2s+1)}}i\biggr]
\end{equation}

\section{Series representations of $\zeta(2s+1)$ and $\beta(2s+1)$} Combining \ref{bet} and \ref{last} we deduce
\begin{equation}\label{ZetaEulerplus}
{(2^{2s+1}-1)\over{2^{2s+1}}}\biggr[\zeta(2s+1)+{{{{ \bf{ E_{2s} }} (\pi i)^{2s+1}}} \over{2(2^{2s+1}-1) \Gamma(2s+1)}}i\biggr] = 2\sum_{k=1}^{\infty} {1\over{(4k-1)^{2s+1}}} 
\end{equation}
and 
\begin{equation}\label{ZetaMinusplus}
{(2^{2s+1}-1)\over{2^{2s+1}}}\biggr[\zeta(2s+1)-{{{{ \bf{ E_{2s} }} (\pi i)^{2s+1}}} \over{2 (2^{2s+1}-1)\Gamma(2s+1)}}i\biggr] = 2\sum_{k=0}^{\infty} {1\over{(4k+1)^{2s+1}}} 
\end{equation} 
because
\begin{equation}
{(2^{2s+1}-1)\over{2^{2s+1}}}\zeta(2s+1) = \sum_{k=0}^{\infty} {1\over{(4k+1)^{2s+1}}}+\sum_{k=1}^{\infty} {1\over{(4k-1)^{2s+1}}}. 
\end{equation} 
We find  
\begin{equation}
\beta(2s+1) =-{(2^{2s+1}-1)\over{2^{2s+1}}}\biggr[{{{{ \bf{ E_{2s} }} (\pi i)^{2s+1}}} \over{2(2^{2s+1}-1)\Gamma(2s+1)}}i\biggr] = \sum_{k=0}^{\infty} {1\over{(4k+1)^{2s+1}}}-\sum_{k=1}^{\infty} {1\over{(4k-1)^{2s+1}}} ;
\end{equation} 
\begin{equation}
{2^{2s+1}\over{2^{2s+1}-1}}\beta(2s+1) =-\biggr[{{{{ \bf{ E_{2s} }} (\pi i)^{2s+1}}} \over{2(2^{2s+1}-1)\Gamma(2s+1)}}i\biggr] ={2^{2s+1}\over{2^{2s+1}-1}}( \sum_{k=0}^{\infty} {1\over{(4k+1)^{2s+1}}}-\sum_{k=1}^{\infty} {1\over{(4k-1)^{2s+1}}}) ;
\end{equation} 
due to the following identities:
\begin{equation}
\zeta(2s+1) = {2^{2s+1}\over{2^{2s+1}-1}}(\sum_{k=0}^{\infty} {1\over{(4k+1)^{2s+1}}}+\sum_{k=1}^{\infty} {1\over{(4k-1)^{2s+1}}});
\end{equation} 
\begin{equation}
\zeta(2s+1)+{2^{2s+1}\over{2^{2s+1}-1}}\beta(2s+1) = {2^{2s+1}\over{2^{2s+1}-1}}\sum_{k=0}^{\infty} {2\over{(4k+1)^{2s+1}}};
\end{equation} 
\begin{equation}
\zeta(2s+1)-{2^{2s+1}\over{2^{2s+1}-1}}\beta(2s+1) = {2^{2s+1}\over{2^{2s+1}-1}}\sum_{k=1}^{\infty} {2\over{(4k-1)^{2s+1}}}.
\end{equation} 

\section{Summary of results}
We confirm the following identities to be valid:
\begin{equation}\label{zetaF}
\zeta(2s) = (-1)^{2s}({{( \psi^{(2s-1)}({1\over4})+\psi^{(2s-1)}({3\over4}))}\over{2^{2s}(2^{2s}-1)}}){1 \over{\Gamma(2s)}} = {{ \pi {d^{(2s-1)}\over{dz^{(2s-1)}}}cot(\pi z) \mid{_{z\rightarrow{1\over4}}}} \over 
                        {{2^{2s}(2^{2s}-1)\Gamma(2s)}} };
\end{equation}
\begin{equation}\label{betaV}
{2^{2s+1}\over{2^{2s+1}-1}}\beta(2s+1) = (-1)^{2s+1}({{( \psi^{(2s)}({1\over4})-\psi^{(2s)}({3\over4}))}\over{2^{2s+1}(2^{2s+1}-1)}}){1 \over{\Gamma(2s+1)}} = {{ \pi {d^{(2s)}\over{dz^{(2s)}}}cot(\pi z) \mid{_{z\rightarrow{1\over4}}}} \over 
                        {{2^{2s+1}(2^{2s+1}-1)\Gamma(2s+1)}} };
\end{equation}

\begin{equation}\label{betaF}
\beta(2s+1) =(-1)^{2s+1}({{( \psi^{(2s)}({1\over4})-\psi^{(2s)}({3\over4}))}\over{2^{2s+1}(2^{2s+1})}}){1 \over{\Gamma(2s+1)}}={{ \pi {d^{(2s)}\over{dz^{(2s)}}}cot(\pi z) \mid{_{z\rightarrow{1\over4}}}} \over 
                        {{2^{2s+1}(2^{2s+1})\Gamma(2s+1)}} };
\end{equation}

\begin{equation}\label{EulerF}
 E_{2s} =(-1)^{2s+1}({{( \psi^{(2s)}({1\over4})-\psi^{(2s)}({3\over4}))}\over{(2\pi i)^{2s+1}}}).2i = {{ \pi {d^{(2s)}\over{dz^{(2s)}}}cot(\pi z) \mid{_{z\rightarrow{1\over4}}}} \over {(2\pi i)^{2s+1}}}.2i;
\end{equation}

\begin{equation}\label{zetaF}
\zeta(2s+1) = (-1)^{2s+1}({{( \psi^{(2s)}({1\over4})+\psi^{(2s)}({3\over4}))}\over{2^{2s+1}(2^{2s+1}-1)}}){1 \over{\Gamma(2s+1)}}
\end{equation}
according to \cite{Ido12b}and \cite{Ido12a}.

\section{Conclusion}
The following identities are derived:
\begin{equation}
\zeta(2s+1) = \biggr[{2^{2s+1}\over{2^{2s+1}-1}}\sum_{k=0}^{\infty} {2\over{(4k+1)^{2s+1}}}\biggr] -\biggr[{{{{ E_{2s}} (\pi i)^{2s+1}}} \over{2(2^{2s+1}-1)\Gamma(2s+1)}}i\biggr];
\end{equation} 
\begin{equation}
\zeta(2s+1)= \biggr[{2^{2s+1}\over{2^{2s+1}-1}}\sum_{k=1}^{\infty} {2\over{(4k-1)^{2s+1}}}\biggr] +\biggr[{{{{  E_{2s}} (\pi i)^{2s+1}}} \over{2(2^{2s+1}-1)\Gamma(2s+1)}}i\biggr]
\end{equation} 
from the general formula and relation
\[
{(2^{2s+1}-1)\over{2^{2s+1}}} \zeta(2s+1)=\beta(2s+1)+{(2^{2s+1}-1)\over{2^{2s+1}}}\biggr[\zeta(2s+1)+{{{{  E_{2s} } (\pi i)^{2s+1}}} \over{2 (2^{2s+1}-1)\Gamma(2s+1)}}i\biggr]
\]
where
\[
\beta(2s+1) =-{(2^{2s+1}-1)\over{2^{2s+1}}}\biggr[{{{{ \bf{ E_{2s} }} (\pi i)^{2s+1}}} \over{2(2^{2s+1}-1)\Gamma(2s+1)}}i\biggr] 
\] 
and
\[
 \zeta(2s+1)={2^{2s+1}\over(2^{2s+1}-1)}\beta(2s+1)+\biggr[{{{{ \bf{ E_{2s} }} (\pi i)^{2s+1}}} \over{{{2(2^{2s+1}-1)}} \Gamma(2s+1)}} - \zeta(2s+1)i\biggr]i
\]
is the closed-form expression for $\zeta(2n+1)$ such that 
\[
{2^{2s+1}\over 2^{2s+1}-1}\sum_{k=1}^{\infty} {2\over{(4k-1)^s}}=
\biggr[{{{{ \bf{ E_{2s} }} (\pi i)^{2s+1}}} \over{{{2(2^{2s+1}-1)}} \Gamma(2s+1)}} - \zeta(2s+1)i\biggr]i
\]

\end{document}